\documentclass[12pt, draft]{article}
\usepackage{amsmath, amsfonts,amsthm, amssymb}
\newtheorem{theorem}{Theorem}[subsection]
\newtheorem{proposition}{Proposition}[subsection]
\newtheorem{corollary}{Corollary}[subsection]
\newtheorem{definition}{Definition}[subsection]
\begin{document}
\title{\bf THE GEOMETRY OF K-ORBITS OF A SUBCLASS OF MD5-GROUPS AND
FOLIATIONS FORMED BY THEIR GENERIC K-ORBITS}
\author{{\bf Le Anh Vu}*\,\, {\bf and}\,\, {\bf Duong Minh Thanh}**\\
** Department of Mathematics and Informatics\\
University of Pedagogy, Ho Chi Minh City, Vietnam\\e-mail:leanhvu@hcmup.edu.vn\\
*** Department of Physics\\
University of Pedagogy, Ho Chi Minh City, Vietnam\\e-mail:
dmthanh@hcmup.edu.vn} \footnotetext{

2000AMS Mathematics Subject Classification: Primary 22E45, Secondary 46E25, 20C20.}
\date{}
\maketitle

\begin{abstract}
The present paper is a continuation of Le Anh Vu's ones [13], [14],
[15]. Specifically, the paper is concerned with the subclass of
connected and simply connected MD5-groups such that their
MD5-algebras $\mathcal{G}$ have the derived ideal ${\mathcal{G}}^{1}
: = [ \mathcal{G},\mathcal{G} ]\equiv$ ${\bf{R}}^{3}$. We shall
describe the geometry of K-orbits of these MD5-groups. The
foliations formed by K-orbits of maximal dimension of these
MD5-groups and their measurability are also presented in the paper.
\end{abstract}
\vskip1cm
\subsection*{Introduction}

    It is well known that the Kirillov's method of orbits (see [3]) plays
the most important in the theory of representations of Lie groups.
By this method, we can obtain all the unitary irreducible
representations of solvable and simply connected Lie groups. The
importance of Kirillov's method of orbits is the co-adjoint
representation (K-representation). Therefore, it is meaningful to
study the K-representation in the theory of representations of
solvable Lie groups.

    The structure of solvable Lie groups and Lie algebras is not  to complicated,
    although the complete classification of them is unresolved up to now. In 1980,
    studying the Kirillov's method of orbits, D. N. Diep (see [2], [6])introduced
    the class of Lie groups and Lie algebras MD. Let G be an n-dimensional Lie group.
    It is called an MDn-group(see [2], [6]), iff its orbits in the co-adjoint
    representation (K-orbits) are orbits of dimension zero or maximal dimension.
    The corresponding Lie algebra are called MDn-algebra. It is worth noticing that,
    for each connected MDn-group, the family of all K-orbits of maximal dimension
    formes a foliation. Therefore, we can combine the study of MDn-algebras and
    MDn-groups with the method of Connes (see [1]) in the theory of foliations and
    operator algebras.

    All MD4-algebras were first listed by D.V Tra in 1984 (see [7]) and then
    classified up to an isomorphism by the first author in 1990 (see [10], [11]).
    The description of the geometry of K-orbits of all indecomposable MD4-groups,
    the topological classification of foliations formed by K-orbits of maximal
    dimension and the characterization of C*-algebras associated to these
    foliations by the method of K-functors were also given by the first author
    in 1990 (see [8], [9], [10], [11]). In 2001, Nguyen Viet Hai (see [3])
    introduced deformation quantization on K-orbits of MD4-groups and obtained
    all representations of MD4-groups. Until now, no complete classification
    of MDn-algebras with $n\ge 5$ is known.

    Recently, the first author continued study MD5-groups and MD5-algebras
$\mathcal{G}$ in cases ${\mathcal{G}}^{1} : = [\mathcal{G},
\mathcal{G}] \equiv$ ${\bf{R}}^{k}$; $k = 1, 2$. (see [13], [14]).
In the present paper we concern with a similar problem for different
MD5-groups and MD5-algebras $\mathcal{G}$ in the case
${\mathcal{G}}^{1}\equiv {\bf{R}}^{3}$. We begin our discussion in
Section 2 by repeating this subclass of MD5-algebras which is listed
in [15] by the first author. Section 3 is devoted to the geometric
description of K-orbits of MD5-groups corresponding to these
MD5-algebras and a discussion of the foliations formed by their
maximal dimensional K-orbits. At first, we recall in Section 1 some
preliminary results and notations which will be used later. For
details we refer the reader to References [1], [2], [4].
\section{Preliminaries}

\subsection{The co-adjoint Representation and K-orbits of a Lie Group}

 Let G be a Lie group. We denote by $\mathcal{G}$ the Lie
algebra of G and by $\mathcal{G^{*}}$ the dual space of
$\mathcal{G}$. To each element g of G we associate an automorphism
$$A_{(g)}:G\longrightarrow G$$
$$\qquad \qquad \qquad\qquad\qquad x\longmapsto A_{(g)}(x):\, =\;gxg^{-1}.$$
(which is called the internal automorphism associated to g).
$A_{(g)}$ induces the tangent map
$${A_{(g)}}_{*}:\mathcal{G}\longrightarrow \mathcal{G}\qquad \qquad \qquad \qquad $$
$$\qquad\qquad \qquad \qquad\qquad X\longmapsto {A_{(g)}}_{*}(X):\, =\;\frac{d}{dt}
[g.exp(tX)g^{-1}]\mid_{t=0}.$$

\begin{definition} The action
$$Ad :G\longrightarrow Aut(\mathcal{G})$$
$$\qquad\qquad\qquad g\longmapsto Ad(g):\, =\;  {A_{(g)}}_{*}$$
is called the adjoint representation of G in $\mathcal{G}$.
\end{definition}
\begin{definition} The action
$$K:G\longrightarrow Aut(\mathcal{G}^{*})$$
$$g\longmapsto K_{(g)}$$
such that
$$\langle K_{(g)}F,X\rangle :\,=\langle F, Ad(g^{-1})X\rangle ;
\quad(F\in {\mathcal{G}}^{*}, X \in \mathcal{G})$$ is called the
co-adjoint representation of G in $\mathcal{{G}^{*}}$.
\end{definition}
\begin{definition} Each orbit of the co-adjoint representation of
G is called a K-orbit. The dimension of a K-orbit of G is always
even. \end{definition}

\subsection{Foliations and Measurable Foliations}

Let V be a smooth manifold. We denote by TV its tangent bundle, so
that for each x $\in$ V, $T_{x}V$ is the tangent space of V at x.
\begin{definition} A smooth subbundle $\mathcal{F}$ of TV is called integrable iff
the following condition is satisfied: every x from V is contained in
a submanifold W of V such that $T_{p}$(W) = ${\mathcal{F}}_{p}$
($\forall p$ $\in$ W).\end{definition}

\begin{definition} A foliation (V, $\, \mathcal {F}$ ) is
given by a smooth manifold V and an integrable subbundle
$\mathcal{F}$ of TV. Then, V is called the foliated manifold and
$\mathcal{F}$ is called the subbundle defining the
foliation.\end{definition}

\begin{definition} The leaves of the foliation (V, $\,
\mathcal{F}$) are the maximal connected submanifolds L of V with
$T_{x}(L) = {\mathcal{F}}_{x}$ ($\forall x \in L $).\end{definition}

    The set of leaves with the quotient topology is denoted by
    V/${\mathcal{F}}$ and called the \emph{space of leaves} of
    (V, $\, \mathcal{F}$). It is a fairly untractable topological space.

   The partition of V in leaves :V = $\bigcup_{\alpha \in V/\mathcal{F}}L_{\alpha}$
   is characterized geometrically by the following local triviality: Every
   $x \in V$ has a system of local coordinates  $\{ U; x^{1}, x^{2}, ..., x^{n} \}
   (x \in U; n = dim \mathcal{F})$  so that for any leaf
   $L$ with $L \cap$ U $\ne \emptyset$ , each connected component of
   $L \cap U$ (which is called a \emph{plaque} of the leaf $L$)
   is given by the equations
$$x^{k+1} = c^{1},\, x^{k+2} = c^{2},\, ...\,, x^{n}\, =\, c^{n-k};\, k = dim{\mathcal{F}}$$
where $c^{1}, c^{2}, ..., c^{n-k}$ are constants (depending on
each plaque). Such a system $\{ U, x^{1}, x^{2}, ..., x^{n} \}$ is called a
\emph{foliation chart}.

   A foliation can be given by a partition of V in a family
$\mathcal{C}$ of its submanifolds such that each $L\in{\mathcal{C}}$
is a maximal connected integral submanifold of some integrable
subbundle $\mathcal{F}$ of TV. Then $\mathcal{C}$ is the family of
leaves of the foliation (V, $\, \mathcal{F}$). Sometimes
$\mathcal{C}$ is identified with $\mathcal{F}$ and we will say that
(V, $\, \mathcal{F}$) is formed by $\mathcal{C}$.

\begin{definition} A submanifold N of the foliated manifold V is
called a transversal iff $T_{x}V = T_{x}N \oplus {\mathcal{F}}_{x}
(\forall x \in N)$. Thus,
 dimN = n - dim$\mathcal{F}$ = codim$\mathcal{F}$.

    A Borel subset B of V such that $B\cap L$ is countable for any leaf L
    is called a Borel transversal to
    (V,\,$\mathcal{F}$).\end{definition}

\begin{definition} A transverse measure $\Lambda$ for the
foliation (V,\,$\mathcal{F}$) is $\sigma$ - additive map B$\mapsto
\Lambda$ (B) from the set of all Borel transversals to [0,
+$\infty$] such that the following conditions are satisfied :

    (i) If $\psi$ : $B_{1} \rightarrow B_{2}$ is a Borel bijection and
    $\psi$(x) is on the leaf of any x$\in B_{1}$, then $\Lambda(B_{1}) = \Lambda(B_{2})$.

    (ii) $\Lambda(K)< + \infty$ if K is any compact subset of a smooth transversal
    submanifold of V.

    By a measurable foliation we mean a foliation
    (V,\,$\mathcal{F}$) equipped with some transverse measure
    $\Lambda$.\end{definition}

 Let (V,\,$\mathcal{F}$) be a foliation with $\mathcal{F}$ is
oriented. Then the complement of zero section of the bundle
${\Lambda}^{k}(\mathcal{F})$ (k = dim$\mathcal{F}$)  has two
components ${\Lambda}^{k}{(\mathcal{F})}^{-}$ and
${\Lambda}^{k}{(\mathcal{F})}^{+}$.

    Let $\mu$ be a measure on V and  $\{ U, x^{1}, x^{2}, ..., x^{n} \}$
    be a foliation chart of (\nobreak V,\,$\mathcal{F}$\nobreak).
    Then U can be identified with the direct product $N \times {\Pi}$
    of some smooth transversal submanifold N of V and a some plaque
$\Pi$. The restriction of $\mu$ on $U \equiv N \times {\Pi}$ becomes the product
${\mu}_{N} \times {\mu}_{\Pi}$ of measures ${\mu}_{N}$ and ${\mu}_{\Pi}$ respectively.

     Let X $\in C^{\infty}{\bigl({\Lambda}^{k}(\mathcal{F})\bigr)}^{+}$
     be a smooth k-vector field and ${\mu}_{X}$ be the measure on each leaf L
     determined by the volume element X.

\begin{definition} The measure $\mu$ is called
X-invariant iff ${\mu}_{X}$ is proportional to ${\mu}_{\Pi}$ for an
arbitrary foliation chart $\{ U, x^{1}, x^{2}, ..., x^{n}
\}$.\end{definition}

 Let (X, $\mu$), (Y, $\nu$) be two pairs where X,Y
 $\in C^{\infty}{\bigl({\Lambda}^{k}(\mathcal{F})\bigr)}^{+}$ and $\mu, \nu$
 are measures on V such that $\mu$ is X-invariant, $\nu$ is Y-invariant.

\begin{definition} ( X, $\mu$ ), ( Y, $\nu$ ) are
equivalent iff Y = $\varphi$ X and $\mu$ = $\varphi \nu$ for some
$\varphi \in C^{\infty}(V).$\end{definition}

    There is one bijective map between the set of transverse measures for
(V,\, $\mathcal{F}$) and the one of equivalence classes of pairs
(X, \, $\mu$), where X $\in
C^{\infty}{\bigl({\Lambda}^{k}(\mathcal{F})\bigr)}^{+}$ and
$\mu$ is a X-invariant measure on V.

    Thus, to prove that (V,\,$\mathcal{F}$) is measurable, we only need choose
    some suitable pair (X, $\mu$) on V.

\section{A Subclass of Indecomposable \\MD5-Algebras and MD5-Groups}
From now on, G will denote a connected simply-connected solvable Lie
group of  dimension 5. The Lie algebra of G is denoted by
$\mathcal{G}$. We always choose a fixed basis $(\nobreak X_{1},
X_{2}, X_{3}, X_{4}, X_{5} \nobreak)$ in $\mathcal{G}$. Then Lie
algebra $\mathcal{G}$ isomorphic to ${\bf R}^{5}$ as a real vector
space. The notation ${\mathcal{G}}^{*}$ will mean the dual space of
$\mathcal{G}$. Clearly ${\mathcal{G}}^{*}$ can be identified with
${\bf R}^{5}$ by fixing in it the basis $(\nobreak X_{1}^{*},
X_{2}^{*}, X_{3}^{*}, X_{4}^{*}, X_{5}^{*}\nobreak )$ dual to the
basis $(\nobreak X_{1}, X_{2}, X_{3}, X_{4}, X_{5} \nobreak)$.

    Recall that a group G is called a {\emph{MD5-group}} iff its K-orbits
    are orbits of dimension zero or maximal dimension. Then its Lie algebra is called
    {\emph{MD5 - algebra}}. Note that for any MDn - algebra
    ${\mathcal{G}}_{0}$ ($ 0 < n <5 $), the direct sum
    $\mathcal{G}$ = ${\mathcal{G}}_{0}$ $\oplus$ ${\bf R}^{5-n}$ of
    ${\mathcal{G}}_{0}$ and the commutative Lie algebra ${\bf R}^{5-n}$ is a MD5-algebra.
    It is called a {\emph{decomposable} } MD5 - algebra, the study of which can be
    directly reduced to the case of MDn - algebras with ($0<n<5$).
    Therefore, we will restrict on the case of {\emph{indecomposable}} MD5 - algebras.

\subsection{List of Indecomposable \\ MD5 - Algebras and MD5 - Groups}

For the sake of convenience, we shall continue using the first
author's notations in [10], [11], [12], [15]. Specifically, we
consider the set
$\lbrace{\mathcal{G}}_{5,3,1(\lambda_1,\lambda_2)}$,
${\mathcal{G}}_{5,3,2(\lambda)}$, ${\mathcal{G}}_{5,3,3(\lambda)}$
${\mathcal{G}}_{5,3,4}$, ${\mathcal{G}}_{5,3,5(\lambda)}$,
${\mathcal{G}}_{5,3,6(\lambda)}$, ${\mathcal{G}}_{5,3,7}$,
${\mathcal{G}}_{5,3,8(\lambda,\varphi)}\rbrace$ of solvable Lie
algebras of dimension 5 which are listed by the first author in
[15]. Each algebra ${\mathcal{G}}$ from this set has
$${\mathcal{G}}^{1} = [ {\mathcal{G}}, {\mathcal{G}} ] =
{\bf R}.X_{3}\oplus {\bf R}.X_{4} \oplus {\bf R}.X_{5}\equiv {\bf R}^{3};\, \,
[X_{1}, X_{2}] = X_{3};\, \, ad_{{X}_1} = 0.$$
The operator $ad_{{X}_{2}} \in End({\mathcal{ G}}^{1})
\equiv  Mat( 3, {\bf R})$ is given  as follows:
\vskip0.5cm
                \begin{description}
                    \item[1.]${\mathcal{G}}_{5,3,1({\lambda}_{1}, {\lambda}_{2})}$ :
                    $$ad_{{X}_2} = \begin{pmatrix} {{\lambda}_1}&0&0\\
                    0&{{\lambda}_2}&0\\0&0&1 \end{pmatrix}; \quad
                     {\lambda}_1, {\lambda}_2 \in {\bf R}\setminus \lbrace
                     1\rbrace, \, {\lambda}_1 \neq {\lambda}_2 \neq 0 .$$ \vskip 0.5cm
                     \item[2.]${\mathcal{G}}_{5,3,2(\lambda)}$ :
                    $$ad_{{X}_2} = \begin{pmatrix} 1&0&0\\
                    0&1&0\\0&0&{\lambda} \end{pmatrix}; \quad
                    {\lambda} \in {\bf R}\setminus \lbrace 0, 1 \rbrace .$$ \vskip 0.5cm
                    \item[3.]${\mathcal{G}}_{5,3,3(\lambda)}$ :
                    $$ad_{{X}_2} = \begin{pmatrix} {\lambda}&0&0\\
                    0&1&0\\0&0&1 \end{pmatrix}; \quad
                    {\lambda} \in {\bf R}\setminus \lbrace 1 \rbrace .$$ \vskip 0.5cm
                    \item[4.]${\mathcal{G}}_{5,3,4}$ :
                    $$ad_{{X}_2} = \begin{pmatrix} 1&0&0\\
                    0&1&0\\0&0&1 \end{pmatrix}.$$ \vskip 0.5cm
                    \item[5.]${\mathcal{G}}_{5,3,5(\lambda)}$ :
                    $$ad_{{X}_2} = \begin{pmatrix} {\lambda}&0&0\\
                    0&1&1\\0&0&1 \end{pmatrix}; \quad
                    {\lambda} \in {\bf R}\setminus \lbrace 1 \rbrace .$$ \vskip 0.5cm
                    \item[6.]${\mathcal{G}}_{5,3,6(\lambda)}$ :
                    $$ad_{{X}_2} = \begin{pmatrix} 1&1&0\\
                    0&1&0\\0&0&{\lambda} \end{pmatrix}; \quad
                    {\lambda} \in {\bf R}\setminus \lbrace 0, 1 \rbrace .$$ \vskip 0.5cm
                    \item[7.]${\mathcal{G}}_{5,3,7}$ :
                    $$ad_{{X}_2} = \begin{pmatrix} 1&1&0\\
                    0&1&1\\0&0&1 \end{pmatrix}.$$ \vskip 0.5cm
                    \item[8.]${\mathcal{G}}_{5,3,8(\lambda, \varphi)}$ :
                    $$ad_{{X}_2} = \begin{pmatrix} cos{\varphi}&-sin{\varphi}&0\\
                    sin{\varphi}&cos{\varphi}&0\\0&0&\lambda \end{pmatrix}; \quad
                     \lambda \in {\bf R}\setminus \lbrace
                     0\rbrace, \, \varphi \in (0, \pi) .$$ \vskip 0.5cm
                \end{description}
\vskip0.5cm
   So we obtain a set of connected and simply-connected solvable Lie groups
   corresponding to the set of Lie algebras listed above. For convenience,
   each such Lie group is also denoted by the same indices as its Lie algebra.
   For example, ${G}_{5,3,6(\lambda)}$ is the connected and simply-connected Lie
group corresponding to ${\mathcal{G}}_{5,3,6(\lambda)}$.
\subsection{Remarks}

In [15], the first author proved that all of the Lie algebras and Lie groups
as were listed above are indecomposable MD5-algebras and MD5-groups. However,
this assertion can be verified by the pictures of K-orbits of considered
Lie groups in the next section.

\section{The Main Results}

\subsection{The Geometry of K-orbits of considered Lie groups}

Throughout this section, G will denote one of the groups
${G}_{5,3,1(\lambda_1,\lambda_2)}$, ${G}_{5,3,2(\lambda)}$,
${G}_{5,3,3(\lambda)}$ ${G}_{5,3,4}$, ${G}_{5,3,5(\lambda)}$,
${G}_{5,3,6(\lambda)}$, ${G}_{5,3,7}$,
${G}_{5,3,8(\lambda,\varphi)}$, ${\mathcal{G}}$ for its Lie algebra,
${{\mathcal{G}} = < X_{1}, X_{2}, X_{3}, X_{4}, X_{5}>}$.
${\mathcal{G}}^{*}$ = $<X_{1}^{*}$, $X_{2}^{*}$, $X_{3}^{*}$,
$X_{4}^{*}$, $X_{5}^{*}> \equiv {\bf R}^{5}$ is the dual space of
${\mathcal{G}}$, $F = {\alpha}X_{1}^{*} + {\beta}X_{2}^{*} +
{\gamma}X_{3}^{*} + {\delta}X_{4}^{*} + {\sigma}X_{5}^{*} \equiv
({\alpha}, {\beta}, {\gamma}, {\delta}, {\sigma})$ an arbitrary
element of ${\mathcal{G}}^{*}$, and finally ${\Omega}_{F}$ the
K-orbit of G which contains $F$. \vskip 0.5cm
\begin{proposition}
    If $G = {G}_{5,3,1(\lambda_1,\lambda_2)}$
    then the picture of the K-orbits of G be described as follows:
        \begin{enumerate}
            \item If\, ${\gamma} = {\delta} = {\sigma} = 0$\, then
$${\Omega}_{F} = \lbrace F( {\alpha}, {\beta}, 0, 0, 0)\rbrace .$$
(K-orbit of dimension zero).
            \item If\, ${\gamma} = {\delta} = 0, {\sigma}\neq 0$\, then
$${\Omega}_{F} = \lbrace F({\alpha}, y, 0, 0, s):  {\sigma}s > 0\rbrace $$
(a part of 2-dimensional plane).
            \item If\, ${\gamma} = 0, {\delta} \neq 0, {\sigma} = 0$\, then
            $${\Omega}_{F} = \lbrace F({\alpha}, y, 0, t, 0):  {\delta}t > 0\rbrace $$
(a part of 2-dimensional plane).
            \item If\, ${\gamma} = 0, {\delta} \neq 0, {\sigma} \neq 0$\, then
$${\Omega}_{F} = \lbrace F({\alpha}, y, 0, t, s):  t =
{\delta}{(\frac{s}{{\sigma}})^{{\lambda}_{2}}} , {\sigma}s > 0\rbrace $$
(a 2-dimensional cylinder).
            \item If\, ${\gamma} \neq 0, {\delta} = {\sigma} = 0$\, then
$${\Omega}_{F} = \lbrace F(x, y, z, 0, 0):
{\lambda}_{1}x = {\lambda}_{1}{\alpha} + {\gamma} - z, {\gamma}z >
0\rbrace $$ (a part of 2-dimensional plane).
            \item If\, ${\gamma} \neq 0, {\delta} = 0, {\sigma} \neq 0$\, then
\begin{eqnarray*}{\Omega}_{F} = \lbrace F(x, y, z, 0, s):
{\lambda}_{1}x = {\lambda}_{1}{\alpha} + {\gamma} - z,
{\lambda}_{1}x = {\lambda}_{1}{\alpha} + {\gamma}(1 -
(\frac{s}{{\sigma}})^{{\lambda}_{1}}),\\ {\sigma}s > 0\rbrace
\end{eqnarray*} (a 2-dimensional cylinder).
            \item If\, ${\gamma} \neq 0, {\delta} \neq 0, {\sigma} = 0$\, then
\begin{eqnarray*}{\Omega}_{F} = \lbrace F(x, y, z, t, 0): {\lambda}_{1}x =
{\lambda}_{1}{\alpha} + {\gamma} - z, {\lambda}_{1}x =
{\lambda}_{1}{\alpha} + {\gamma}({1 -
(\frac{t}{{\delta}})^{\frac{{\lambda}_{1}}{{\lambda}_{2}}}}),\\
{\delta}t > 0\rbrace \end{eqnarray*} (a 2-dimensional cylinder).
           \item If\, ${\gamma} \neq 0, {\delta} \neq 0, {\sigma} \neq 0$\, then
\begin{eqnarray*}{\Omega}_{F} = \lbrace F(x, y, z, t, s):
{\lambda}_{1}x = {\lambda}_{1}{\alpha} + {\gamma} - z,
{\lambda}_{1}x = {\lambda}_{1}{\alpha} + {\gamma}
({1 - (\frac{s}{{\sigma}})^{{\lambda}_{1}}}), \\
t = {\delta}(\frac{s}{{\sigma}})^{{\lambda}_{2}}, {\sigma}s >
0\rbrace \end{eqnarray*} (a 2-dimensional cylinder).
        \end{enumerate}
\end{proposition}
{\bf Sketch of the proof of Propositions 3.1.1}

   For G, we denote the set $\lbrace F_{U} \in {\mathcal{G}}^{*}
/ U \in {\mathcal{G}}\rbrace$ by  ${\Omega}_{F}(\mathcal{G})$, where $F_{U}$
is the linear form on the Lie algebra $\mathcal{G}$ of G defined by
$$\langle F_{U}, A \rangle = \langle F, exp({ad}_{U})(A)\rangle, A, U \in \mathcal{G}.$$
Let $U = a.X_{1} + b.X_{2} + c.X_{3} + d.X_{4} + f.X_{5}$ be an arbitrary of $\mathcal{G}$;
where $ a, b, c, d, f \in {\bf R}$.
Upon direct computation, we get :

$$exp({ad}_{U})= \begin{pmatrix}
  1 & 0 & 0 & 0 & 0 \\
  0 & 1 & 0 & 0 & 0 \\
  -{\sum_{n=1}^{\infty}\frac{b^{n}{\lambda_{1}}^{n-1}}{n!}}&
  (a-c\lambda_{1})\sum_{n=1}^{\infty}\frac{(b{\lambda_{1}})^{n-1}}{n!}
   & e^{b{\lambda_{1}}} & 0 & 0 \\
  0 & -d\sum_{n=1}^{\infty}\frac{b^{n-1}{\lambda_{1}}^{n}}{n!} & 0 & e^{b\lambda_{2}}& 0 \\
  0 & -f\sum_{n=1}^{\infty}\frac{b^{n-1}}{n!}& 0 & 0& e^{b} \\\end{pmatrix}.$$
\vskip0.3cm

Thus, $F_{U}$ is given as follows:
\vskip 0.3cm
   $x = \alpha - \gamma{\sum_{n=1}^{\infty}\frac{b^{n}{\lambda_{1}}^{n-1}}{n!}}; $
\vskip 0.3cm

   $y = \beta + \gamma (a-c\lambda_{1})\sum_{n=1}^{\infty}
   \frac{(b{\lambda_{1}})^{n-1}}{n!} - \delta d\sum_{n=1}^{\infty}
   \frac{b^{n-1}{\lambda_{1}}^{n}}{n!} - \sigma f\sum_{n=1}^{\infty}\frac{b^{n-1}}{n!}; $
\vskip 0.3cm

   $z = \gamma e^{b\lambda_{1}};$
\vskip 0.3cm

   $t = \delta e^{b\lambda_{2}};$
\vskip 0.3cm

   $s = \sigma e^{b}.$
\vskip0.3cm

So ${\Omega}_{F}(\mathcal{G})$ is described and the equation
${\Omega}_{F}(\mathcal{G}) = {\Omega}_{F}$ is verified by the same
method presented in [8], [10], and [12]. \hfill{$\square$}

   According to the method of the proof of Proposition 1, we get the following results.

\begin{proposition}

    If $G = {G}_{5,3,2({\lambda})}$ then the picture of the K-orbits of G be
    described as follows:
       \begin{enumerate}
            \item If\, ${\gamma} = {\delta} = {\sigma} = 0$\, then
$${\Omega}_{F} = \lbrace F( {\alpha}, {\beta}, 0, 0, 0)\rbrace.$$
(K-orbit of dimension zero).
            \item If\, ${\gamma} = {\delta} = 0, {\sigma}\neq 0$\, then
$${\Omega}_{F} = \lbrace F({\alpha}, y, 0, 0, s):  {\sigma}s > 0\rbrace $$
(a part of 2-dimensional plane).
            \item If\, ${\gamma} = 0, {\delta} \neq 0, {\sigma} = 0$\, then
$${\Omega}_{F} = \lbrace F({\alpha}, y, 0, t, 0):  {\delta}t > 0\rbrace $$
(a part of 2-dimensional plane).
            \item If\, ${\gamma} = 0, {\delta} \neq 0, {\sigma} \neq 0$\, then
$${\Omega}_{F} = \lbrace F({\alpha}, y, 0, t, s):
s = {\sigma}(\frac{t}{{\delta}})^{\lambda} , {\delta}t > 0\rbrace $$
(a 2-dimensional cylinder).
            \item If\, ${\gamma} \neq 0, {\delta} = {\sigma} = 0$\, then
$${\Omega}_{F} = \lbrace F(x, y, z, 0, 0):  x = {\alpha} + {\gamma} - z,
{\gamma}z > 0\rbrace $$
(a part of 2-dimensional plane).
            \item If\, ${\gamma} \neq 0, {\delta} = 0, {\sigma} \neq 0$\, then
$${\Omega}_{F} = \lbrace F(x, y, z, 0, s):  x = {\alpha} + {\gamma} - z,
s = {\sigma}(\frac{z}{{\gamma}})^{\lambda}, {\gamma}z > 0\rbrace $$
(a 2-dimensional cylinder).
            \item If\, ${\gamma} \neq 0, {\delta} \neq 0, {\sigma} = 0$\, then
$${\Omega}_{F} = \lbrace F(x, y, z, t, 0):  x = {\alpha} + {\gamma} - z,
x = {\alpha} + (1 - \frac{t}{{\delta}}){\gamma}, {\delta}t > 0\rbrace $$
(a part of 2-dimensional plane).
            \item If\, ${\gamma} \neq 0, {\delta} \neq 0, {\sigma} \neq 0$\, then
$${\Omega}_{F} = \lbrace F(x, y, z, t, s):
x = {\alpha} + {\gamma} - z, x = {\alpha} + (1 -
\frac{t}{{\delta}}){\gamma}, s =
{\sigma}(\frac{t}{{\delta}})^{\lambda}, {\delta}t > 0\rbrace $$ (a
2-dimensional cylinder). \hfill{$\square$}
        \end{enumerate}
\end{proposition}

\begin{proposition}

    If $G = {G}_{5,3,3({\lambda})}$ then the picture of the K-orbits of G
    be described as follows:
        \begin{enumerate}
            \item If\, ${\gamma} = {\delta} = {\sigma} = 0$\, then
$${\Omega}_{F} = \lbrace F( {\alpha}, {\beta}, 0, 0, 0)\rbrace.$$
(K-orbit of dimension zero).
            \item If\, ${\gamma} = {\delta} = 0, {\sigma}\neq 0$\, then
$${\Omega}_{F} = \lbrace F({\alpha}, y, 0, 0, s):  {\sigma}s > 0\rbrace $$
(a part of 2-dimensional plane).
            \item If\, ${\gamma} = 0, {\delta} \neq 0, {\sigma} = 0$\, then
$${\Omega}_{F} = \lbrace F({\alpha}, y, 0, t, 0):  {\delta}t > 0\rbrace $$
(a part of 2-dimensional plane).
            \item If\, ${\gamma} = 0, {\delta} \neq 0, {\sigma} \neq 0$\, then
$${\Omega}_{F} = \lbrace F({\alpha}, y, 0, t, s):
{\delta}s = {\sigma}t , {\delta}t > 0\rbrace $$
(a 2-dimensional cylinder).
            \item If\, ${\gamma} \neq 0, {\delta} = {\sigma} = 0$\, then
$${\Omega}_{F} = \lbrace F(x, y, z, 0, 0):  {\lambda}x = {\lambda}{\alpha} + {\gamma} - z,
{\gamma}z > 0\rbrace $$
(a part of 2-dimensional plane).
            \item If\, ${\gamma} \neq 0, {\delta} = 0, {\sigma} \neq 0$\, then
$${\Omega}_{F} = \lbrace F(x, y, z, 0, s):
{\lambda}x = {\lambda}{\alpha} + {\gamma} - z,
{\lambda}x = {\lambda}{\alpha} + {\gamma}(1 - (\frac{s}{{\sigma}})^{\lambda}),
{\sigma}s > 0\rbrace $$
(a 2-dimensional cylinder).
            \item If\, ${\gamma} \neq 0, {\delta} \neq 0, {\sigma} = 0$\, then
$${\Omega}_{F} = \lbrace F(x, y, z, t, 0):
{\lambda}x = {\lambda}{\alpha} + {\gamma} - z,
z = {\gamma}(\frac{t}{{\delta}})^{{\lambda}}, {\delta}t > 0\rbrace $$
(a 2-dimensional cylinder).
            \item If\, ${\gamma} \neq 0, {\delta} \neq 0, {\sigma} \neq 0$\, then
\begin{eqnarray*}{\Omega}_{F} = \lbrace F(x, y, z, t, s):
{\lambda}x = {\lambda}{\alpha} + {\gamma} - z, {\lambda}x =
{\lambda}{\alpha} + {\gamma}(1 - (\frac{t}{{\delta}})^{\lambda}),\\
{\sigma}t = {\delta}s, {\delta}t > 0\rbrace \end{eqnarray*} (a
2-dimensional cylinder). \hfill{$\square$}
        \end{enumerate}
\end{proposition}

\begin{proposition}

     If $G = {G}_{5,3,4}$ then the picture of the K-orbits of
     G be described as follows:
        \begin{enumerate}
            \item If\, ${\gamma} = {\delta} = {\sigma} = 0$\, then
$${\Omega}_{F} = \lbrace F( {\alpha}, {\beta}, 0, 0, 0)\rbrace.$$
(K-orbit of dimension zero).
            \item If\, ${\gamma} = {\delta} = 0, {\sigma}\neq 0$\, then
$${\Omega}_{F} = \lbrace F({\alpha}, y, 0, 0, s):  {\sigma}s > 0\rbrace $$
(a part of 2-dimensional plane).
            \item If\, ${\gamma} = 0, {\delta} \neq 0, {\sigma} = 0$\, then
$${\Omega}_{F} = \lbrace F({\alpha}, y, 0, t, 0):  {\delta}t > 0\rbrace $$
(a part of 2-dimensional plane).
            \item If\, ${\gamma} = 0, {\delta} \neq 0, {\sigma} \neq 0$\, then
$${\Omega}_{F} = \lbrace F({\alpha}, y, 0, t, s):  {\delta}s ={\sigma}t ,
{\delta}t > 0\rbrace $$
(a part of 2-dimensional plane).
            \item If\, ${\gamma} \neq 0, {\delta} = {\sigma} = 0$\, then
$${\Omega}_{F} = \lbrace F(x, y, z, 0, 0):  x = {\alpha} + {\gamma} - z,
{\gamma}z > 0\rbrace $$
(a part of 2-dimensional plane).
            \item If\, ${\gamma} \neq 0, {\delta} = 0, {\sigma} \neq 0$\, then
$${\Omega}_{F} = \lbrace F(x, y, z, 0, s):  x = {\alpha} + {\gamma} - z, x =
{\alpha} + {\gamma}(1 - \frac{s}{{\sigma}}), {\sigma}s > 0\rbrace $$
(a part of 2-dimensional plane).
            \item If\, ${\gamma} \neq 0, {\delta} \neq 0, {\sigma} = 0$\, then
$${\Omega}_{F} = \lbrace F(x, y, z, t, 0):  x = {\alpha} + {\gamma} - z,
z = {\gamma}\frac{t}{{\delta}}, {\delta}t > 0\rbrace $$
(a part of 2-dimensional plane).
            \item If\, ${\gamma} \neq 0, {\delta} \neq 0, {\sigma} \neq 0$\, then
$${\Omega}_{F} = \lbrace F(x, y, z, t, s):  x = {\alpha} + {\gamma} - z, x = {\alpha} +
{\gamma}(1 - \frac{s}{{\sigma}}), {\sigma}t = {\delta}s, {\delta}t >
0\rbrace $$ (a part of 2-dimensional plane). \hfill{$\square$}
        \end{enumerate}
\end{proposition}

\begin{proposition}

     If $G = {G}_{5,3,5({\lambda})}$ then the picture of the
     K-orbits of G be described as follows:
        \begin{enumerate}
            \item If\, ${\gamma} = {\delta} = {\sigma} = 0$\, then
$${\Omega}_{F} = \lbrace F( {\alpha}, {\beta}, 0, 0, 0)\rbrace.$$
(K-orbit of dimension zero).
            \item If\, ${\gamma} = {\delta} = 0, {\sigma}\neq 0$\, then
$${\Omega}_{F} = \lbrace F({\alpha}, y, 0, 0, s):  {\sigma}s > 0\rbrace $$
(a part of 2-dimensional plane).
            \item If\, ${\gamma} = 0, {\delta} \neq 0, {\sigma} = 0$\, then
$${\Omega}_{F} = \lbrace F({\alpha}, y, 0, t, s):  s = tln\frac{t}{{\delta}},
{\delta}t > 0\rbrace $$
(a 2-dimensional cylinder).
            \item If\, ${\gamma} = 0, {\delta} \neq 0, {\sigma} \neq 0$\, then
$${\Omega}_{F} = \lbrace F({\alpha}, y, 0, t, s):
s = {\sigma}\frac{t}{{\delta}} + tln\frac{t}{{\delta}}, {\delta}t > 0\rbrace $$
(a 2-dimensional cylinder).
            \item If\, ${\gamma} \neq 0, {\delta} = {\sigma} = 0$\, then
$${\Omega}_{F} = \lbrace F(x, y, z, 0, 0):
{\lambda}x = {\lambda}{\alpha} + {\gamma} - z,
{\gamma}z > 0\rbrace $$
(a part of 2-dimensional plane).
            \item If\, ${\gamma} \neq 0, {\delta} = 0, {\sigma} \neq 0$\, then
$${\Omega}_{F} = \lbrace F(x, y, z, 0, s):
{\lambda}x = {\lambda}{\alpha} + {\gamma} - z,
{\lambda}x = {\lambda}{\alpha} + {\gamma}(1 - (\frac{s}{{\sigma}})^{\lambda}),
{\sigma}s > 0\rbrace $$ (a
2-dimensional cylinder).
            \item If\, ${\gamma} \neq 0, {\delta} \neq 0, {\sigma} = 0$\, then
$${\Omega}_{F} = \lbrace F(x, y, z, t, s):
{\lambda}x = {\lambda}{\alpha} + {\gamma} - z, z = {\gamma}(\frac{t}{{\delta}})^{{\lambda}},
s = tln\frac{t}{{\delta}}, {\delta}t > 0\rbrace $$(a 2-dimensional cylinder).
            \item If\, ${\gamma} \neq 0, {\delta} \neq 0, {\sigma} \neq 0$\, then
\begin{eqnarray*}{\Omega}_{F} = \lbrace F(x, y, z, t, s):
{\lambda}x = {\lambda}{\alpha} + {\gamma} - z,
{\lambda}y = {\lambda}{\alpha} + {\gamma}(1 - (\frac{t}{{\delta}})^{\lambda}),\\
s = {\sigma}\frac{t}{{\delta}} + tln\frac{t}{{\delta}}, {\delta}t >
0\rbrace \end{eqnarray*} (a 2-dimensional cylinder). \hfill{$\square$}
        \end{enumerate}
\end{proposition}

\begin{proposition}

     If $G = {G}_{5,3,6({\lambda})}$ then the picture of the K-orbits of G be
described as follows:
        \begin{enumerate}
            \item If\, ${\gamma} = {\delta} = {\sigma} = 0$\, then
$${\Omega}_{F} = \lbrace F( {\alpha}, {\beta}, 0, 0, 0)\rbrace.$$
(K-orbit of dimension zero).
            \item If\, ${\gamma} = {\delta} = 0, {\sigma}\neq 0$\, then
$${\Omega}_{F} = \lbrace F({\alpha}, y, 0, 0, s):  {\sigma}s > 0\rbrace $$
(a part of 2-dimensional plane).
            \item If\, ${\gamma} = 0, {\delta} \neq 0, {\sigma} = 0$\, then
$${\Omega}_{F} = \lbrace F({\alpha}, y, 0, t, 0):  {\delta}t > 0\rbrace $$
(a part of 2-dimensional plane).
            \item If\, ${\gamma} = 0, {\delta} \neq 0, {\sigma} \neq 0$\, then
$${\Omega}_{F} = \lbrace F({\alpha}, y, 0, t, s):
s = {\sigma}(\frac{t}{{\delta}})^{{\lambda}}, {\delta}t > 0\rbrace $$
(a 2-dimensional cylinder).
            \item If\, ${\gamma} \neq 0, {\delta} = {\sigma} = 0$\, then
$${\Omega}_{F} = \lbrace F(x, y, z, t, 0):
x = {\alpha} + {\gamma} - z, t = zln\frac{z}{{\gamma}}, {\gamma}z > 0\rbrace $$
(a 2-dimensional cylinder).
            \item If\, ${\gamma} \neq 0, {\delta} = 0, {\sigma} \neq 0$\, then
$${\Omega}_{F} = \lbrace F(x, y, z, t, s):  x = {\alpha} + {\gamma} - z, t =
zln\frac{z}{{\gamma}}, s = {\sigma}(\frac{z}{{\gamma}})^{{\lambda}},
{\sigma}s > 0\rbrace $$ (a 2-dimensional cylinder).
            \item If\, ${\gamma} \neq 0, {\delta} \neq 0, {\sigma} = 0$\, then
$${\Omega}_{F} = \lbrace F(x, y, z, t, 0):  x = {\alpha} + {\gamma} - z,
t = zln\frac{z}{{\gamma}} + {\delta}\frac{z}{{\gamma}}, {\gamma}z > 0\rbrace $$
(a 2-dimensional cylinder).
            \item If\, ${\gamma} \neq 0, {\delta} \neq 0, {\sigma} \neq 0$\, then
$${\Omega}_{F} = \lbrace F(x, y, z, t, s):  x = {\alpha} + {\gamma} - z,
t = zln\frac{z}{{\gamma}} + {\delta}\frac{z}{{\gamma}}, s =
{\sigma}(\frac{z}{{\gamma}})^{{\lambda}}, {\gamma}z > 0\rbrace $$ (a
2-dimensional cylinder). \hfill{$\square$}
        \end{enumerate}
\end{proposition}

\begin{proposition}

     If $G = {G}_{(5,3,7)}$ then the picture of the K-orbits of G
     be described as follows:
        \begin{enumerate}
            \item If\, ${\gamma} = {\delta} = {\sigma} = 0$\, then
$${\Omega}_{F} = \lbrace F( {\alpha}, {\beta}, 0, 0, 0)\rbrace.$$
(K-orbit of dimension zero).
            \item If\, ${\gamma} = {\delta} = 0, {\sigma}\neq 0$\, then
$${\Omega}_{F} = \lbrace F({\alpha}, y, 0, 0, s):  {\sigma}s > 0\rbrace $$
(a part of 2-dimensional plane).
            \item If\, ${\gamma} = 0, {\delta} \neq 0, {\sigma} = 0$\, then
$${\Omega}_{F} = \lbrace F({\alpha}, y, 0, t, s):  s = tln\frac{t}{{\delta}},
{\delta}t > 0\rbrace $$
(a 2-dimensional cylinder).
            \item If\, ${\gamma} = 0, {\delta} \neq 0, {\sigma} \neq 0$\, then
$${\Omega}_{F} = \lbrace F({\alpha}, y, 0, t, s):
s = tln\frac{t}{{\delta}} + {\sigma}\frac{t}{{\delta}}, {\delta}t > 0\rbrace $$
(a 2-dimensional cylinder).
            \item If\, ${\gamma} \neq 0, {\delta} = {\sigma} = 0$\, then
$${\Omega}_{F} = \lbrace F(x, y, z, t, s):  x = {\alpha} + {\gamma} - z,
t = zln\frac{z}{{\gamma}}, s = \frac{z}{2}ln^{2}\frac{z}{{\gamma}},{\gamma}z > 0\rbrace $$ (a 2-dimensional
cylinder).
            \item If\, ${\gamma} \neq 0, {\delta} = 0, {\sigma} \neq 0$\, then
$${\Omega}_{F} = \lbrace F(x, y, z, t, s):
x = {\alpha} + {\gamma} - z, t = zln\frac{z}{{\gamma}},
s = \frac{z}{2}ln^{2}\frac{z}{{\gamma}} + {\sigma}\frac{z}{{\gamma}},
{\gamma}z > 0\rbrace $$
(a 2-dimensional cylinder).
            \item If\, ${\gamma} \neq 0, {\delta} \neq 0, {\sigma} = 0$\, then
\begin{eqnarray*}{\Omega}_{F} = \lbrace F(x, y, z, t, s):  x = {\alpha} + {\gamma}
- z, t = zln\frac{z}{{\gamma}} + {\delta}\frac{z}{{\gamma}}, \\
s = \frac{z}{2}ln^{2}\frac{z}{{\gamma}} + \frac{{\delta}}{{\gamma}}zln\frac{z}{{\gamma}},
{\gamma}z > 0\rbrace \end{eqnarray*}
(a 2-dimensional cylinder).
            \item If\, ${\gamma} \neq 0, {\delta} \neq 0, {\sigma} \neq 0$\, then
\begin{eqnarray*}{\Omega}_{F} = \lbrace F(x, y, z, t, s):  x = {\alpha} + {\gamma}
- z, t = zln\frac{z}{{\gamma}} + {\delta}\frac{z}{{\gamma}}, \\
s = \frac{z}{2}ln^{2}\frac{z}{{\gamma}} + \frac{{\delta}}{{\gamma}}zln\frac{z}{{\gamma}}
 + \sigma\frac{z}{{\gamma}}, {\gamma}z > 0\rbrace \end{eqnarray*}
(a 2-dimensional cylinder). \hfill{$\square$}
           \end{enumerate}
\end{proposition}

\begin{proposition}

    If $G = {G}_{(5,3,8({\lambda},{\varphi}))}$ then the picture of the
K-orbits of G be described as follows:
        \begin{enumerate}
            \item If\, ${\gamma} = {\delta} = {\sigma} = 0$\, then
$${\Omega}_{F} = \lbrace F( {\alpha}, {\beta}, 0, 0, 0)\rbrace.$$
(K-orbit of dimension zero).
            \item If\, ${\gamma} = {\delta} = 0, {\sigma}\neq 0$\, then
            $${\Omega}_{F} = \lbrace F({\alpha}, y, 0, 0, s):  {\sigma}s > 0\rbrace $$
(a part of 2-dimensional plane).
            \item If\, ${\gamma}^{2} + {\delta}^{2} \neq 0 $\, then
$${\Omega}_{F} = \lbrace F(x, y, z + it, s) = F(x, y, {\gamma}e^{b{e}^{-i{\varphi}}}
+ {\delta}e^{b{e}^{i{\varphi}}}, {\sigma}e^{{\lambda}b}) \rbrace $$
(a 2-dimensional cylinder). \hfill{$\square$}
        \end{enumerate}
\end{proposition}

\subsection{Remark and Corollary}

      Note that if $G$ is not ${G}_{5,3,8(\lambda,\varphi)}$ then $G$
      is exponential (see [5]). Hence ${\Omega}_{F}= {\Omega}_{F}(\mathcal{G})$.

   For ${G}_{5,3,8(\lambda,\varphi)}$, the equation
${\Omega}_{F}= {\Omega}_{F}(\mathcal{G})$ is verified by using [12, Lemma II.1.5 ].

     As an immediate consequence of above propositions we have the following corollary.

\begin{corollary}
   All\,  of \, ${G}_{5,3,1(\lambda_1,\lambda_2)}$,\,  ${G}_{5,3,2(\lambda)}$,\,
 ${G}_{5,3,3(\lambda)}$, \,   ${G}_{5,3,4}$,\,  ${G}_{5,3,5(\lambda)}$,\,
 ${G}_{5,3,6(\lambda)}$,\,  ${G}_{5,3,7}$,\,  ${G}_{5,3,8(\lambda,\varphi)}$
 are MD5-groups. \hfill{$\square$} \end{corollary}

\subsection{MD5-Foliations Associated to \\Considered MD5-Groups}

\begin{theorem}
   Let G $\in \lbrace$ {${G}_{5,3,1(\lambda_1,\lambda_2)}$, ${G}_{5,3,2(\lambda)}$,
${G}_{5,3,3(\lambda)}$, ${G}_{5,3,4}$, ${G}_{5,3,5(\lambda)}$,
${G}_{5,3,6(\lambda)}$, ${G}_{5,3,7}$,
${G}_{5,3,8(\lambda,\varphi)}$}$\rbrace$, ${\mathcal{F}}_{G}$ be the
family of all its K-orbits of maximal dimension and $V_{G} =
\bigcup{\lbrace \Omega / \Omega \in {\mathcal{F}}_{G}\rbrace}$ .
Then $( V_{G},\,{\mathcal{F}}_{G} )$ is a measurable foliation in
the sense of Connes. We call it MD5-foliation associated to G.
\end{theorem} {\bf Sketch of the Proof of Theorem 3.3.1}

The proof is analogous to the case of MD4-groups in [8], [10], [12]
or the first examples of MD5-group in [13], [14]. First, we need to
define a smooth tangent 2-vector field on the manifold $V_{G}$ such
that each K-orbit $\Omega$ from ${\mathcal{F}}_{G}$ is a maximal
connected integrable submanifold corresponding to it. As the next
step, we have to show that the Lebegues measure is invariant for
that 2-vector field. This steps can be complete by direct
computations. We therefore omit detail comutations of the proof.
\hfill{$\square$}

\subsection{Concluding Remark}

We conclude this paper with the following comments.

   1. It should be note that the results of Propositions 3.1.1 - 3.1.8
   and Theorem 3.3.1 still hold for all indecomposable connected
   (no-simply connected) MD5-groups correponding to above-mentioned MD5-algebras.

   2. There are some interesting questions can be raised for further study.
   That are giving the topological classification of considered MD5-foliations,
   characterizing $C^{*}$-algebras associated to these foliations, constructing
   deformation quantization on K-orbits of considered MD5-groups to obtain all
   representations of them, ...  The authors will announce the new results
   of this questions in the next papers.

\subsection*{Acknowledgement}

    The authors would like take this opportunity to thank their teacher,
    Prof. DSc. Do Ngoc Diep for his excellent advice and support.
    Thanks are due also to the first author's colleague, Prof. Nguyen Van Sanh
    for his encouragement.

\end{document}